\documentclass[12pt]{amsart}
\usepackage{amsmath,amsthm,amsfonts,amssymb}

\usepackage{eucal}
\newcommand{\1}{\mathbf{1}}

\DeclareMathOperator{\SL}{SL}
\DeclareMathOperator{\Hom}{Hom}
\DeclareMathOperator{\im}{i}

\newcommand{\Th}{\theta}
\newcommand{\tS}{{\tilde S}}
\newcommand{\ts}{{\tilde s}}

\newcommand{\one}{\mathbf{1}}

\newcommand{\C}{\mathbb C}
\newcommand{\Z}{\mathbb Z}\newcommand{\R}{\mathbb R}
\newcommand{\Q}{\mathbb Q}
\newcommand{\CC}{\mathcal{C}}

\newcommand{\Gal}{\mathrm{Gal}}
\newcommand{\N}{\mathbb N}
\newcommand{\ot}{\otimes}

\newcommand{\lan}{\langle}
\newcommand{\ra}{\rangle}
\newcommand{\ep}{\epsilon}

\numberwithin{equation}{section}
\newtheorem{thm}[equation]{Theorem}
\newtheorem{theorem}[equation]{Theorem}

\newtheorem{lem}[equation]{Lemma}

\newtheorem{conj}[equation]{Conjecture}

\newtheorem{prop}[equation]{Proposition}
\theoremstyle{definition}

\newtheorem{remark}[equation]{Remark}

\newtheorem{case}[equation]{Case}

\begin{document}
\title[Non-self-dual modular categories]
{On the classification of non-self-dual modular categories}
\begin{abstract}
We develop a symbolic computational approach to classifying low-rank modular categories.
 We use this technique to classify pseudo-unitary modular categories of rank at most $5$ that are non-self-dual, i.e. those for which some object is not isomorphic to its dual object.
\end{abstract}
\author{Seung-moon Hong}
\email{seungmoon.hong@utoledo.edu}
\address{Department of Mathematics\\
     University of Toledo\\
    Toledo, OH 43606\\
    U.S.A.}

\author{Eric Rowell}
\email{rowell@math.tamu.edu}
\address{Department of Mathematics\\
    Texas A\&M University \\
    College Station, TX 77843\\
    U.S.A.}
\thanks{The second author is partially supported by NSA grant H98230-08-1-0020.
We thank Victor Ostrik, Zhenghan Wang, Jennifer Franko and Matthew Young for useful comments.}

\maketitle
\section{Introduction}
This paper is concerned with the classification of modular categories as defined by Turaev (see \cite{BK}).  This problem has been considered in different guises going back (at least) to the early 1990s.  Some early results are found in the physics literature where classifications have been obtained for certain restricted classes of \textit{modular fusion algebras}: the finite rank unital based ring \cite{O1} describing the tensor product (fusion) rules of a modular category.  Gepner and Kapustin \cite{GK} determined those modular fusion algebras of rank $\leq 6$ with very small structure constants (e.g. $\leq 1$ for rank $6$).  Fuchs \cite{Fuchs} classified those of rank $\leq 3$ under certain (physically inspired) compatibility conditions, while Eholzer \cite{Eh} classified modular fusion algebras of rank $\leq 4$ that are \emph{strongly-modular}, i.e. the kernel of the modular representation 
is assumed to contain a congruence subgroup.  A very recent result of Ng and Schauenburg \cite{NS} shows that this assumption is superfluous. 

Classifying up to fusion algebras leaves only finite ambiguity due to a generalized form of Ocneanu rigidity due to Etingof, Nikshych and Ostrik (see \cite[Remark 2.33]{ENO}):
\begin{prop}
 There are at most finitely many braided fusion categories realizing a given fusion algebra.
\end{prop}
It is well known that the number of compatible ribbon structures for a given braided fusion category is finite (see e.g. \cite{Kas}), so that the result holds for modular categories as well.

A recent conjecture due to Wang is:
\begin{conj}\cite[Conjecture 6.1]{RSW}\label{wangconj}
There are finitely many modular categories of rank $n$ for any fixed $n\in\N$.
\end{conj}
In \cite{RSW} a classification of modular categories of rank $n\leq 4$ is given under the assumption of pseudo-unitarity (which verifies Conjecture \ref{wangconj} for these cases).  In fact, the cases categories of rank $n=2,3$ have been classified in greater generality by Ostrik \cite{O2,O3}.  Moreover, the conjecture has been verified (for fusion categories) in \cite[Proposition 8.38]{ENO} under the assumption that the global dimension is integral.  

For fixed rank $n$, the classification breaks naturally into two types: 1) those categories for which every object is isomorphic to its dual, and 2) the categories for which at least one object is not isomorphic to its dual.
The main results of this paper provide a classification of pseudo-unitary modular categories of rank $n\leq 5$ of type 2).  We choose to focus on the non-self-dual case for several reasons, besides the fact that there are fewer parameters than in the self-dual case.  Firstly, one of the main results used in \cite{RSW} for the self-dual, rank $n\leq 4$ case (\cite[Theorem 2.14]{RSW}) does not apply to the non-self-dual case.  Secondly, we recover the results of \cite[Appendix A]{RSW} for ranks $n\leq 4$ with significantly less effort, employing Proposition \ref{rswprop} to greater advantage as well as a new general result (Theorem \ref{maxsdtheorem}). Finally, our approach uses symbolic computation techniques that can be generalized, whereas \cite{RSW} relies upon (difficult \emph{ad hoc}) hand computations.  

The problem of classifying modular categories can be viewed in parallel with that of classifying finite simple groups, and is an interesting theoretical problem from this perspective.
Our original motivation for classifying low-rank modular categories comes from topological quantum computation \cite{FLW}.  Modular categories play a key role in the mathematical description of topological phases of matter--the physical systems upon which topological quantum computation are expected to be built.  See \cite[Section 6]{RSW} for details on this relationship.  From a symbolic computation perspective, classifying rank $n$ modular fusion algebras corresponds to solving a large system of polynomials in three sets of variables:  $O(n^2)$ algebraic integers, $O(n^3)$ non-negative integers and $n$ roots of unity (of undetermined order).  For $n>2$ the system resists a purely machine computation solution, both because of the system's complexity and the essential way in which integrality and Galois groups are used.

Here is a more detailed description of the paper.  In Section \ref{prelim} we obtain some general results on integral modular categories and apply them to the classification of rank $n=3,4$ non-self-dual categories.  Section \ref{rank5} gives the classification of rank $n=5$ non-self-dual modular fusion algebras.  In Section \ref{conclusions} we summarize the technique and discuss its potential applications and generalizations.

\section{Preliminaries}\label{prelim}
A modular category is a finite semisimple $\C$-linear rigid braided balanced monoidal category with simple identity object satisfying a certain non-degeneracy condition (invertibility of the $S$-matrix).  For our purposes the precise axioms satisfied by a modular category will not be important; we refer the interested reader to the text \cite{BK} for details.  We shall mostly be interested in the polynomial consequences of the axioms, which we will describe below.  \textbf{We assume throughout that $\dim(X_i)>0$ and refer to this condition as \textit{pseudo-unitarity}}.  This differs slightly from the terminology in \cite{ENO}, but conforms with the assumptions in \cite{RSW}.  

Let $\CC$ be a rank $n$ modular category with inequivalent simple objects represented by classes $\{\one=X_0,X_1,\ldots,X_{n-1}\}$.  We will abuse notation and refer to $X_i$ as a simple object as this will cause no confusion. The (class of) the object dual to $X_i$ will be denoted either $X^*_i$ or $X_{i^*}$. We denote by $Gr(\CC)$ the Grothendieck ring of $\CC$, which is a finite rank unital based ring in the notation of \cite{O1} and more specifically is a \emph{modular fusion algebra}.  Since $\CC$ is braided, $Gr(\CC)$ must be commutative.

Let $\tS$ denote the $S$-matrix of $\CC$ with entries
$\ts_{ij}$.  Notice that we have $\ts_{0i}=\dim(X_i)$, and $\tS^2=D^2C$ where $C_{ij}=\delta_{i,j^*}$ and $D^2=\dim(\CC)=\sum_i \dim(X_i)^2$.  Define $\psi_i$ to be the character of $Gr(\CC)$ corresponding to the normalized $i$th column of $\tS$, i.e. $\psi_i(X_j)=\ts_{ij}/\dim(X_i)$.  Observe that these characters are orthogonal since the $\tS$-matrix is (projectively) unitary, i.e.
$$\sum_k \psi_i(X_k)\overline{\psi_j(X_k)}=\delta_{ij}\dim(\CC)/\dim(X_i)^2.$$
We let $K=\Q(\ts_{ij})$ and define $\Gal(\CC):=\Gal(K/\Q)$. 

We will make extensive use of the following results from \cite[Theorem 2.10]{RSW} originally due to Coste and Gannon \cite{CG}:
\begin{prop}\label{rswprop}
 Let $Gal(\CC)$ be the Galois group of a rank $n$ modular category $\CC$.
Then
\begin{enumerate}
\item  The action of $\Gal(\CC)$ on $K$ induces an injective group homomorphism $Gal(\CC)\rightarrow S_n$ acting by permutations of $\psi_i$.
\item Let $\sigma$ be the image of any element of $Gal(\CC)$ under the above homomorphism.  Then there exists $\ep_{i,\sigma}=\pm 1$ such that
\begin{equation*}
\tilde{s}_{j,k}=\ep_{\sigma(j),\sigma}\ep_{k,\sigma}\tilde{s}_{\sigma(j),\sigma^{-1}(k)}.
\end{equation*}
\item $Gal(\CC)$ is abelian.
\end{enumerate}

\end{prop}

Let $N_i$ denote the fusion matrix corresponding to the simple object $X_i$, that is
$N_i$ has $(k,j)$-entry $N_{i,j}^k=\dim\Hom(X_i\ot X_j,X_k)$, so that $X_i\rightarrow N_i$ induces the left-regular representation of $Gr(\CC)$. We point out that this differs from the notation of \cite{RSW} in that the fusion matrices $N_i$ are transposed.  The $\{N_i\}$ is a \emph{commutative} set of $n\times n$ $\N$-matrices since $Gr(\CC)$ is commutative. There is an important relationship between $\tS$ and the $N_i$ known as the \emph{Verlinde formula}, which implies that the fusion coefficients $N_{i,j}^k$ are rational functions in the entries of $\tS$.  This is a consequence of the fact that the $\tS$ matrix simultaneously diagonalizes the family of matrices $\{N_i\}$.  The eigenvalues of $N_i$ are $\{\psi_j(X_i): 0\leq j\leq n-1\}$ so that the entries of $\psi_i$ are algebraic integers (as the characteristic polynomial of $N_i$ is monic with integer coefficients).  From this we see that if $\psi_i$ is $\Gal(\CC)$-fixed then $\ts_{ij}/\dim(X_i)\in\Z$.  In particular if $\psi_0$ is fixed then $\CC$ is \textit{integral}, that is, $\dim(X_i)$ is an integer for each $i$.

One has the following consequence:
\begin{theorem}\label{maxsdtheorem}
 Suppose $\CC$ is a modular category of odd rank $n$ such that the
only self-dual object is $\one$.  Then $\CC$ is integral.
\end{theorem}
\begin{proof}
 Since complex conjugation interchanges the characters $\psi_i$ and $\psi_i^*$ it follows that
$(1\; 2)(3\; 4)\cdots(n-2\; n-1)\in G\cong\Gal(\CC)$.  Since $G$ is abelian, it follows that $G$ must fix $\psi_0$ hence all dimensions of simple objects are integral.
\end{proof}

In the case of integral modular categories we have:
\begin{lem}\label{egyptlemma}
 Suppose $\CC$ is an integral modular category of rank $n$.  
Let $p_1\geq p_2\geq \cdots\geq p_{n}=1$ be the sequence of dimensions of simple objects ordered in a weakly decreasing fashion.  Then:
\begin{enumerate}
 \item[(a)] the numbers $x_i:=\dim(\CC)/(p_i)^2$ form a weakly increasing sequence of integers such that $\sum_{i=1}^n 1/x_i=1$ and
\item[(b)] $i\leq x_i\leq (n-i+1)u_{i}$ where $u_1:=1$ and $u_{k+1}:=u_k(u_k+1)$.
\end{enumerate}
\end{lem}
\begin{proof}
By \cite[Lemma 1.2]{EG} $\dim(\CC)/\dim(X_i)^2$ is an algebraic integer for any modular category $\CC$ so in an integral modular category it must be an ordinary integer.  Since $\dim(\CC):=\sum_{i=1}^n p_i^2$, we obtain (a) by dividing by $\dim(\CC)$.  The lower bound in (b) is clear: if some $x_i<i$ then $1/x_j>1/i$ for each $j\leq i$ so that $\sum_{j=1}^i 1/x_j>1$, contradiction.  On the other hand, it is a classical result \cite{land} that if $\sum_{i=1}^k 1/y_k=r$
where $y_1\leq y_2\leq\cdots\leq y_k$ are integers then $y_i\leq (k-i+1)/r_{i-1}$ where $r_0=r$ and $r_{i+1}:=r_i-1/y_{i+1}$.  So (b) will follow once we have shown that $1/r_{i-1}\leq u_i$ for the special case $r_0=1$.  This is also a classical problem solved in \cite{Cur} where it is shown that $1/u_{m+1}\leq 1-\sum_{i=1}^m 1/y_i$ for any $m$ integers $y_i$.
\end{proof}

\begin{remark}\label{algorithm}
 The shifted sequence $u_i+1$ is known as Sylvester's sequence, the first few terms of which are $2,3,7,43,1807,\ldots$.  Lemma \ref{egyptlemma} suggests an algorithm for finding the possible dimensions of rank $n$ integral modular categories $\CC$: we look for integers $x_1\leq\cdots\leq x_n=\dim(\CC)$ such that
\begin{enumerate}
 \item $\sqrt{x_n/x_i}\in\N$ for each $i$
\item $\sum_i 1/x_i=1$ and
\item the $x_i$ satisfy the inequalities of Lemma \ref{egyptlemma}(b).
\end{enumerate}
 Unfortunately the $u_i$ increase at a doubly-exponential rate ($u_6$ is about 3 million), so this algorithm is quite slow for large $n$.
\end{remark}

The term \emph{modular} is applied to these categories because the axioms require that the matrix $\tS$ is invertible and, moreover, there exists a diagonal matrix $T$ whose diagonal entries $[1,\theta_1,\ldots,\theta_{n-1}]$
are roots of unity so that 

$$\begin{pmatrix}0 & -1\\1 &
0\end{pmatrix}\rightarrow \tS, \quad \begin{pmatrix}1 & 1\\0 &
1\end{pmatrix}\rightarrow T$$ induces a (projective) representation of the modular group $\SL(2,\Z)$.  We shall not use this fact directly, but we will use the useful relation (see \cite[Chapter 3]{BK}):
\begin{equation}\label{threls}
 \theta_i\theta_j\ts_{ij}=\sum_k N_{i^*,j}^{k}\dim(X_k)\theta_k
\end{equation}

\subsection{Ranks 3 and 4}
Rank 3 and 4 non-self-dual modular categories have already been classified in \cite{RSW}.  We include these cases both for completeness and because our proofs are somewhat shorter.
Recall that a fusion category is \emph{pointed} if every simple object is invertible.  In the pseudo-unitary setting this is equivalent to $\dim(X_i)=1$ for any simple object $X_i$.  Pointed modular categories of rank $n$ are easily classified up to fusion rules: their Grothendieck rings are isomorphic to the group algebra of an abelian group of order $n$.

\begin{thm}
 Any rank 3 or 4 non-self-dual modular category is pointed.
\end{thm}
\begin{proof}
Firstly, Theorem \ref{maxsdtheorem} implies that any rank 3 non-self-dual modular category must be integral.   Using Remark \ref{algorithm} we obtain only one solution $x_1=x_2=x_3=3$ which corresponds to $\dim(X_i)=1$ for all $i$--the pointed case.

For the rank 4 case we have the following $\tS$-matrix:
$$\begin{pmatrix} 1 & d & g&g\\
   d & x & a & a\\
g& a & h & \overline{h}\\
g & a & \overline{h} & h
  \end{pmatrix}$$

From the
 formula
 $\tilde{s}_{j,k}=\epsilon_{\tau(j),\tau}\epsilon_{k,\tau}\tilde{s}_{\tau(j),\tau^{-1}(k)}$
 for each $\tau\in \Gal(\CC)$ in Proposition \ref{rswprop}, we find that the only solution (respecting orthogonality of the columns $\tS$)
is $x=1$, $a=-g$.  Using orthogonality of the columns of $\tS$ we immediately obtain $d=g^2$ and $D=g^2+1$.  Since  $X_2\cong X_3^*$ and $X_3$, we find that
$\dim(X_2)\dim(X_3)=g^2=1+m_1g^2+m_2g+m_3g$ where $m_i=N_{2,3}^i$ is the multiplicity of $X_i$ in $X_2\ot X_3$.  Clearly $m_1=0$ and $m_2=m_3$ since $X_2\ot X_3$ is self-dual.  We also have $\dim(X_1)\dim(X_2)=g^3=n_1g^2+n_2g+n_3g$, where $n_i=N_{1,2}^i$ is the multiplicity of $X_i$ in $X_1\ot X_2$.  Now $n_2=m_1=0$ since 
$$0=m_1=N_{2,3}^1=N_{1,2}^2=0$$
as $X_2^*=X_3$.  These facts imply
$g^2=1+2m_2g=n_1g+n_3$.  If $g$ is not integral then $m_2>0$ and any linear relation in $g$ over $\Z$ must vanish identically, so that $(1-n_3)-(2m_2-n_1)g=0$ implies $n_3=1$ and $n_1=2m_2$. Next consider $\dim(X_2)^2=g^2=k_1g^2+(k_2+k_3)g$ where
$k_i=N_{2,2}^i$.  We see that $k_1=n_3=1$, and $k_2=m_3=m_2$ using the symmetries of the $N_{i,j}^k$.  Thus we have $g^2=g^2+(m_2+k_3)g$ hence $m_2=m_3=0$.  But then we have $n_1=2m_2=0$ 
which cannot happen unless $g$ is integral.
Since  $D/g=(g^2+1)/g$ is an algebraic integer by \cite[Lemma 2.1]{EG}, $1/g$ is an integer so $g=1$.  Thus we arrive at the pointed case for rank 4 as well.
\end{proof}

\section{Rank 5}\label{rank5}
Now let us assume that $\CC$ is a modular category of rank $5$.  A classification of all modular fusion algebras of rank $N=5$ with bound fusion multiplicities $N_{i,j}^k\leq 3$ is found in \cite{GK}.  There are two such with non-self-dual objects.  In the Kac-Moody algebra (or rational conformal field theory) formulation these are realized as $SU(5)_1$ and $SU(3)_4/\Z_3$.  
Quantum group realizations can be obtained as semisimple sub-quotients of $Rep(U_q(\mathfrak{sl}_5))$ at $q=e^{\pi \im/6}$ and $Rep(U_q(\mathfrak{sl}_3))$ at $q=e^{\pi \im/7}$ where in the latter case one takes the subcategory generated by objects labeled by integer weights (see \cite{Rsurvey} for details).
The fusion rules for $SU(5)_1$ are the same as the addition in $\Z_5$, while the fusion rules for $SU(3)_4/\Z_3$ are given by:
\begin{eqnarray*}
 &N_1&=\begin{pmatrix}
0&1&0&0&0\\
1&2&1&1&1\\
0&1&1&1&1\\
0&1&1&1&0\\
0&1&1&0&1\\
\end{pmatrix},
N_2=\begin{pmatrix}
0&0&1&0&0\\
0&1&1&1&1\\
1&1&1&0&0\\
0&1&0&0&1\\
0&1&0&1&0\\
\end{pmatrix}\\
&N_3&=\begin{pmatrix}
0&0&0&0&1\\
0&1&1&0&1\\
0&1&0&1&0\\
1&1&0&0&0\\
0&0&1&1&0\\
\end{pmatrix}
\end{eqnarray*}
and $N_4=N_3^T$.  We will adopt Gepner and Kapustin's notation and denote the corresponding modular categories as $SU(5)_1$ and $SU(3)_4/\Z_3$.
We will show that this classification is complete without the assumption $N_{i,j}^k\leq 3$.

We first consider the integral case.
\begin{theorem}\label{ptdtheorem}
 Suppose $\CC$ is a rank $5$ integral modular category.  Then $\CC$ is pointed, i.e. $Gr(\CC)$ is isomorphic to $Gr(SU(5)_1)\cong \Z_5$.
\end{theorem}
\begin{proof} From Remark \ref{algorithm} 
we obtain only two possible solutions: 
$$x_1=x_2=x_3=x_4=x_5=5\quad\text{and} \quad x_1=2, x_2=x_3=x_4=x_5=8.$$  The former corresponds to $\dim(X_i)=1$ for all $i$, and hence is pointed.  The latter would be a modular category with four simple objects of dimension $1$ and one simple object of dimension $2$.  This can be ruled out in the following way: the objects $X_0,\ldots,X_3$ of dimension $1$ must be invertible and form a group of order $4$ so the eigenvalues of $N_i$ $0\leq i\leq 3$ are $4$th roots of unity.  This implies that $\psi_4(X_i)=\ts_{i,4}/2=\pm 1$ for $0\leq i\leq 3$, hence $\ts_{i,4}=\pm 2$ for $0\leq i\leq 3$.  But $\sum_i|\ts_{i,4}|^2=8$ which is a contradiction.  So the only possibility is that $\CC$ is pointed.
\end{proof}

\begin{lem}
Suppose that $\CC$ is a rank $5$ pseudo-unitary modular category with $X_3^*\cong X_4$.  Then either $\CC$ is pointed or $\Gal(\CC)$ is is isomorphic to the one of the two $S_5$-subgroups $\lan (0\; 1\; 2),(3\; 4)\ra$ or $\lan (0\; 1), (3\; 4)\ra$.
\end{lem}

\begin{proof}
 If $\psi_0$ is $\Gal(\CC)$-fixed then $\CC$ is integral, hence pointed by Theorem \ref{ptdtheorem}.  By Theorem \ref{ptdtheorem} any rank $5$ modular category with two dual pairs of simple objects must be integral.  We conclude that $(3\; 4)\in\Gal(\CC)$ since complex conjugation fixes $\psi_i$ with $i\leq 2$ and interchanges $\psi_3$ and $\psi_4$.  Thus we may assume that $\Gal(\CC)$ is an abelian subgroup of $Cent_{S_5}((3\; 4))$ containing $(3\; 4)$ that does not fix $0$.  After relabeling $X_1$ and $X_2$ if necessary, we see that the only possibilities are those given in the statement.
\end{proof}

In either case we may assume that the $\tS$-matrix for a rank 5 pseudo-unitary modular category with one dual pair is:
\begin{equation}\label{s-mat}
 \tS=\begin{pmatrix}
          1 & d & f & g & g\\
d & x & y& a & a\\
f & y & z& b & b\\
g & a & b& h & \overline{h}\\
g & a & b& \overline{h} & h
         \end{pmatrix}
\end{equation}
where $\{d,f,g,x,y,z,a,b\}\subset\R$, $h=h_1+\im h_2\in\C\setminus\R$ and each of $d,f$ and $g$ are at least $1$ and algebraic integers.  We will denote $D^2=\dim(\CC)=1+d^2+f^2+2g^2$.   

Using the symmetries
$$N_{i,j}^k=N_{j,i}^k=N_{i,k^*}^{j^*}=N_{i^*,j^*}^{k^*}$$
where the involution on labels ${}^*$ fixes $1$ and $2$ and $3^*=4$, we obtain
the following fusion matrices (note that $N_4=N_3^{T}$):

\begin{eqnarray}
&N_1&=\begin{pmatrix}
0&1&0&0&0\\
1&n_1&n_2&n_3&n_3\\
0&n_2&n_4&n_5&n_5\\
0&n_3&n_5&n_6&n_7\\
0&n_3&n_5&n_7&n_6\\
\end{pmatrix},\label{N1}\\
&N_2&=\begin{pmatrix}
0&0&1&0&0\\
0&n_2&n_4&n_5&n_5\\
1&n_4&n_8&n_9&n_9\\
0&n_5&n_9&n_{10}&n_{11}\\
0&n_5&n_9&n_{11}&n_{10}\\
\end{pmatrix}\label{N2}\\
&N_3&=\begin{pmatrix}
0&0&0&0&1\\
0&n_3&n_5&n_7&n_6\\
0&n_5&n_9&n_{11}&n_{10}\\
1&n_6&n_{10}&n_{12}&n_{13}\\
0&n_7&n_{11}&n_{14}&n_{12}\\
\end{pmatrix}\label{N3}
\end{eqnarray}

We find that there is a single class of modular categories with $Gal(\CC)=\lan (0\; 1\; 2),(3\; 4)\ra\cong \Z_6$ (up to fusion rules).  Precisely, we have:
\begin{theorem}\label{Z6theorem}
Suppose $\CC$ has $Gal(\CC)=\lan (0\; 1\; 2),(3\; 4)\ra$ with $X_3\cong X_4^*$ and $X_1$ and $X_2$ self-dual.  Then $Gr(\CC)$ is isomorphic to $Gr(SU(3)_4/\Z_3)$.
\end{theorem}

The proof is in several steps, which we outline here:

\begin{enumerate}
 \item[Step 1] Use the Galois group and $\tS^2=D^2C$ to determine $\tS$ in terms of the simple dimensions $d,f,g$.
\item[Step 2] Use the commutivity of the $N_i$ and certain symmetries of their characteristic polynomials gleaned from Step 1 to obtain a single Diophantine equation in 3 non-negative integer variables $\alpha,\beta$ and $\gamma$, and expressions for all $N_{i,j}^k$ as polynomials in $\alpha,\beta$ and $\gamma$.  This uses symbolic computation techniques (Gr\"obner bases).
\item[Step 3] Use eqn. (\ref{threls}) to find a single (reducible) polynomial of degree 3 satisfied by $\theta_2$ over the field $\Q(d)$.  Conclude that $\theta_1$ and $\theta_2$ must be (primitive) roots of unity of degree $1,2,3,4,6,7,9,14$ or $18$.
\item[Step 4] Consider each case from Step 3 to conclude that the Diophantine equation has no solutions unless $\theta_2$ is a $7$th root of unity, in which case it has a unique solution. 
\end{enumerate}

Consider the $\tS$-matrix in eqn. (\ref{s-mat}) and suppose that $Gal(\CC)=\lan (0\; 1\; 2),(3\; 4)\ra$. Let
$\sigma=(0\;1\;2)\in S_5$. From the formula
 $\tilde{s}_{j,k}=\epsilon_{\tau(j),\tau}\epsilon_{k,\tau}\tilde{s}_{\tau(j),\tau^{-1}(k)}$
 for each $\tau\in \Gal(\CC)$ in Proposition \ref{rswprop}, we have
 the following:
\begin{eqnarray*}
&&1=\tilde{s}_{0,0}=\epsilon_{1,\sigma}\epsilon_{0,\sigma}y,\quad
d=\tilde{s}_{1,0}=\epsilon_{2,\sigma}\epsilon_{0,\sigma}z\\
&&g=\tilde{s}_{3,0}=\epsilon_{3,\sigma}\epsilon_{0,\sigma}b,\quad
x=\tilde{s}_{1,1}=\epsilon_{2,\sigma}\epsilon_{1,\sigma}f\\
&&a=\tilde{s}_{3,1}=\epsilon_{3,\sigma}\epsilon_{1,\sigma}g,\quad
a=\tilde{s}_{4,1}=\epsilon_{4,\sigma}\epsilon_{1,\sigma}g\\
&&b=\tilde{s}_{3,2}=\epsilon_{3,\sigma}\epsilon_{2,\sigma}a
\end{eqnarray*}

From these we have four possibilities:

\begin{enumerate}
\item[Case 1] $\epsilon_{0,\sigma}=\epsilon_{1,\sigma}=\epsilon_{2,\sigma}=\epsilon_{3,\sigma}=\epsilon_{4,\sigma}$

\item[Case 2]
$\epsilon_{0,\sigma}=\epsilon_{1,\sigma}\neq\epsilon_{2,\sigma}=\epsilon_{3,\sigma}=\epsilon_{4,\sigma}$

\item [Case 3]
$\epsilon_{0,\sigma}=\epsilon_{2,\sigma}\neq\epsilon_{1,\sigma}=\epsilon_{3,\sigma}=\epsilon_{4,\sigma}$

\item [Case 4]
$\epsilon_{0,\sigma}=\epsilon_{3,\sigma}=\epsilon_{4,\sigma}\neq\epsilon_{1,\sigma}=\epsilon_{2,\sigma}$
\end{enumerate}

In Case 1, we have $y=1,z=d,x=f,x=b=g$. Then $\tS^2=D^2C$ gives
us an easy contradiction that $d+df+f+2g^2=0$.

For Case 2 we have $y=1,z=-d,x=-f,a=b=-g$. Then from $\tS^2=D^2C$
we have $d-df+f-2g^2=0$, which implies that $d+f=df+2g^2\geq df+2$
and thus $-1\geq(d-1)(f-1)\geq 0$, a contradiction.

Cases 3 and 4 are equivalent up to permutation of objects
$X_1$ and $X_2$, so we will focus on Case 3.  We obtain $y=-1,z=d,a=g,b=-g,x=-f$. The
equation $\tS^2=D^2C$ implies the following:

\begin{eqnarray*}
&&1+d^2+f^2+2g^2=D^2\\
&&d-df-f+2g^2=0\\
&&g+dg-fg+2gh_1=0\\
&&3g^2+2(h_1^2+h_2^2)=D^2\\
&&3g^2+2h_1^2-2h_2^2=0
\end{eqnarray*}

Without loss of generality, we may assume $h_2>0$ by conjugating if necessary.  These equations then reduce to 
\begin{eqnarray*}
&&h_2=D/2, h_1=(f-1-d)/2\\
&&D^2=1+d^2+f^2+2g^2\\
&&d-df-f+2g^2=0
\end{eqnarray*}
which describes a variety of dimension 2.  The important facts we obtain are: (a) all of the entries of $\tS$ are determined by $d,f$ and $g$, which (b) satisfy the important relation:
\begin{equation}\label{orthog}
 d-df-f+2g^2=0
\end{equation}

Next we consider the fusion matrices and the Diophantine equations their entries satisfy.
Let us denote by $p_i(z)$ the characteristic polynomial of $N_i$.  From the $\tS$-matrix we have:

\begin{eqnarray*}
 &&p_1(z)=(z-1)^2(z^3+c_1z^2+c_2z-1),\\
 &&p_2(z)=(z+1)^2(z^3+c_2z^2-c_1z+1),\\
&&p_3(z)=(z^2+a_1z+a_2)(z^3+b_1z^2+b_2z+b_3)
\end{eqnarray*}
where $c_1=(f/d-d+1/f)$, $c_2=(1/d-f-d/f)$, $a_1=-2h_1/g$, $a_2=(|h|^2/g)^2$, $b_1=g(1/f-1/d-1)$, $b_2=g^2(1/d-1/df-1/f)$ and $b_3=g^3/df$.

Observe that each of the cubic factors of $p_i$ and $(z^2+a_1z+a_2)$ are irreducible over $\Q$, since we assume that the Galois group cyclically permutes the three roots of each of these cubics, and the roots of the quadratic factor of $p_3$ are not real.  In particular, none of $f,d$ or $g$ is an integer.  Moreover, we have $\Q(d)=\Q(d,f,g)$ and $[\Q(d):\Q]=3$. The second statement is a consequence of $\Q(d)/\Q$ being an abelian Galois extension, and the first statement follows from the facts that $-f/d$ is a Galois conjugate of $d$ and 
$b_1\in\Z$.

We obtain several useful relations among the $n_i$ by computing the $p_1$ and $p_2$ directly from $N_1$ and $N_2$ and comparing coefficients:
\begin{enumerate}
 \item The $z^4$ coefficient of $p_1$ is equal to $-1$ times the $z$ coefficient of $p_2$.
\item The constant terms of $p_1$ and $p_2$ are $-1$ and $1$ respectively.
\item the $z^3$ coefficient of $p_1$ is equal to the $z^2$ coefficient of $p_2$ and the $z^2$ coefficient of $p_1$ is equal to $-1$ times the $z^3$ coefficient of $p_2$.
\item The only linear factors of $p_1$ and $p_2$ are $(z-1)$ and $(z+1)$ respectively. 
\end{enumerate}
For example, the last observation implies that $n_6=n_7+1$ and $n_{11}=n_{10}+1$.
The resulting relations among the $n_i$ together with those implied by the pairwise commutivity of the $N_i$ yield two useful results (using Maple's Gr\"obner basis algorithm):
\begin{enumerate}
 \item All entries can be unique expressed as polynomials in $t:=n_{10},v:=n_{12}$ and $u:=n_{14}$ and
\item $u,v$ and $t$ satisfy the Diophantine equation
\begin{equation}\label{dio}
 (v+u)((2t-1)u-(2t+3)v)+4t^2+2t+1=0
\end{equation}
\end{enumerate}
Explicitly in terms of $u,v$ and $t$ the fusion coefficients are:
\begin{eqnarray*}
&n_1&= \frac{2t(2+(u-v)^2)+(u-v)(u^3-u^2v+3u-v^2u+5v+v^3)}{2},\\
&n_2&=u^2-v^2+2t,\\
&n_3&=\frac{((u+v)(u-v)^2+3v+u)+2t(u-v)}{2},\\
&n_4&=2t+1,\quad n_5=u+v,\\
&n_6&=1/2(u^2-v^2+2t+1),\quad n_7=1/2(u^2-v^2+2t-1),\\
&n_8&=t(2+(u-v)^2)+2v^2-2uv+1,\\
&n_9&=t(u-v)+2v,\quad n_{10}=t,\quad n_{11}=t+1,\\
&n_{12}&=v,\quad n_{13}=v,\quad n_{14}=u
\end{eqnarray*}

Observe that eqn. (\ref{dio}) is an indefinite binary quadratic form for fixed $t$ with discriminant $4(2t+1)^2$ (a square) and for any $t$ there are at most finitely many solutions.  It appears that there are infinitely many integers $t\geq 0$ for which (\ref{dio}) has solutions in $u$ and $v$.  For $t=0$ we have the unique (non-negative) solution $u=1$, $v=0$.  We will eventually show that this is the only solution that is realized by a modular category.

Next we describe the relations among the integers $u,v$ and $t$ and the algebraic integers $d,f$ and $g$.  These come from the fact that $\psi_0=[1,d,f,g,g]$ is a simultaneous eigenvector for the $N_i$ with eigenvalue $\dim(X_i)$.  That is:
\begin{equation}\label{chars}
 N_1\psi_0=d\psi_0,\quad N_2\psi_0=f\psi_0,\quad N_3\psi_0=g\psi_0
\end{equation}

The variety described by the vanishing of the ideal generated by
eqns. (\ref{orthog}),(\ref{dio}) and (\ref{chars}) is 2-dimensional, so these are not sufficient to determine all of $u,v,t,d,f$ and $g$ up to finitely many choices. 

Next we analyze eqns. (\ref{threls}).  Since $X_3^*\cong X_4$ we have that $\theta_4=\theta_3$, so that eqns. (\ref{threls}) introduce 3 new variables $\theta_1,\theta_2$ and $\theta_3$, which are roots of unity. We obtain the following six equations from (\ref{threls}):
\begin{eqnarray}
-f\Th_1^2&=&1+n_1d\Th_1+n_2f\Th_2+2n_3g\Th_3\label{threl1} \\
-\Th_1\Th_2&=&n_2d\Th_1+n_4f\Th_2+2n_5g\Th_3\label{threl2}\\
g\Th_1\Th_3&=&n_3d\Th_1+n_5f\Th_2+(n_6+n_7)g\Th_3\label{threl3}\\
-g\Th_2\Th_3&=&n_5d\Th_1+n_{9}f\Th_2+(n_{10}+n_{11})g\Th_3\label{threl4}\\
d\Th_2^2&=&1+n_4d\Th_1+n_8f\Th_2+2n_9g\Th_3\label{threl5}\\
(f-d-1)\Th_3^2&=&1+(n_6+n_7)d\Th_1+(n_{11}+n_{10})f\Th_2+\label{threl6}\\ 
&&(2n_{12}+n_{13}+n_{14})g\Th_3\nonumber
\end{eqnarray}
The last equation comes from $\ts_{3,3}+\ts_{3,4}=2h_1=(f-d-1)$.

Before proceeding, we note that taking eqns. (\ref{orthog}),(\ref{chars}) and (\ref{threl1}-\ref{threl6}) together with the assumptions $t=0,u=0$ and $v=1$ imply the following:

$$\Th_2^6+\Th_2^5+\Th_2^4+\Th_2^3+\Th_2^2+\Th_2+1=0,\quad \Th_1=\Th_2^3,\quad \Th_3=\overline{\Th_2}$$
and 
$$d^3-3d^2-4d-1=0,\quad f^3-4f^2+3f+1=0,\quad g^3-2g^2-g+1=0.$$

Since $t=0$ and eqn. (\ref{dio}) imply $u=1$ and $v=0$ we see that if we show that the relations (\ref{dio}), (\ref{orthog}), (\ref{chars}) and (\ref{threl1}-\ref{threl6}) together with assumptions that the $\Th_i$ are roots of unity and $u,v$ and $t$ are integral force $t=0$, then Theorem \ref{Z6theorem} will follow.

With this in mind, we first wish to bound the degree of the root of unity $\Th_2$.  We claim that $\Th_2$ satisfies a polynomial of degree at most 2 with coefficients in $\Q(d)$.  

We use (\ref{threl2}) to solve
$$\Th_1=(n_4f\Th_2+2n_5g\Th_3)/(\Th_2+n_2d)$$
 which is possible since 
$\Th_2+n_2d=\Th_2+(u^2-v^2+2t)d\not=0$ since $|d|\not=1$ and $|\Th_2|\not=0$.  We then eliminate $\Th_1$ from the remaining 5 relations and take their numerators.  Under this substitution (\ref{threl4}) remains linear in $\Th_3$ so that we may solve 
$$\Th_3=-\frac{f\Th_2(n_9\Th_2+n_9n_2d-n_5dn_4)}
{g(\Th_2^2+(n_{11}+n_{10}+n_2d)\Th_2-2n_5^2d+n_{11}n_2d+n_{10}n_2d)}.$$ 
If the denominator of this expression vanishes then $\Th_2$ satisfies a monic degree two polynomial with coefficients in $\Q(d)$ and the claim follows.

Otherwise, we may eliminate $\Th_3$ from the remaining 4 relations to obtain  polynomials in $\Q(d)[\Th_2]$.  By factoring these polynomials and removing spurious factors (such as $\Th_2+n_2d$) we get two polynomials of degree 4 (from (\ref{threl3}) and (\ref{threl5})) and two polynomials of degree 5 (from (\ref{threl1}) and (\ref{threl6})).  Expressing the coefficients of the two degree 4 polynomials in terms of $d,f,g,t,u$ and $v$ and reducing modulo the ideal generated by eqns. (\ref{dio}), (\ref{orthog}) and (\ref{chars})
these become $(u+v)\Th_2^4+\cdots+d^2(u+v)$
and $d\Th_2^4+\cdots+d$ where the omitted terms are too long to include here.
We cancel the $\Th_2^4$ terms to obtain a polynomial of degree 3 in $\Th_2$ with coefficients in $\Q(d)$ with constant term $d(d^2-1)(u+v)\not=0$ (since $u$ and $v$ cannot both be zero as otherwise (\ref{dio}) has no integer solutions for $t$).

Since $\Th_2$ satisfies a (non-zero) polynomial of degree at most $3$ over $\Q(d)$, we see that $m:=[\Q(\Th_2,d):\Q(d)]\in\{1,2,3\}$, which is a Galois extention since $\Th_2$ is a root of unity.   We can immediately eliminate $m=3$ since the degree of the non-trivial cyclotomic extension $\ell:=[\Q(\Th_2):\Q]$ cannot divide $9$.  Thus the claim follows so that $\ell\in\{1,2,6\}$ hence
$\Th_2$ must be a primitive root of unity of degree $s\in\{1,2,3,4,6,7,9,14,18\}$
Moreover, since $\Th_3\in\Q(\Th_2,d)$ and $\Th_1\in\Q(\Th_3,\Th_2,d)$ the same is true of $\Th_1$ and $\Th_3$.  

We can now proceed to Step 4.  We systematically eliminate all possibilities except $s=7$.  The strategy is as follows: 
\begin{enumerate}
 \item[(a)] Let $\mathcal{I}$ be ideal generated by eqns. (\ref{dio}), (\ref{orthog}), (\ref{chars}), (\ref{threl1}-\ref{threl6}) and $\varphi_s(\Th_2)$--the cyclotomic polynomial whose roots are primitive $s$th roots of unity.
\item[(b)] Compute a Gr\"obner basis for $\mathcal{I}$ with respect to the lexicographic order induced by $d<f<g<\Th_1<\Th_2<\Th_3<u<v<t$.
\item[(c)] Factor the first polynomial in the output of step (b), which will be a polynomial in $v,t$ or just $t$.
\item[(d)] Conclude that no integer solutions exist (by parity arguments, for example).
\end{enumerate}
We consider each case in turn.
\begin{case}
 Assume that $s=1$, (i.e. $\Th_2=1$).  Then we find that  at least \textit{one of} the following three polynomials must vanish:
\begin{eqnarray*}
 &&(2t+1)(6t+1),\label{eq11}\\
&&16(1+t)^2v^4+(4(2t-1)(1+t)(t^2-4t-2))v^2+(4t^2+2t+1)(t^2-t+1)^2,\label{eq12}\\
&&64v^4+(104+32t^2-128t)v^2+(7-10t+4t^2)^2\label{eq13}\\
\end{eqnarray*}
Clearly the first expression cannot be zero, while the latter two are odd for
$t\in\Z$ hence have no solutions.
\end{case}
\begin{case}
 Assume that $s=2$.  Then the Gr\"obner basis computation gives us the consequence: $(2t+1)(6t+1)(4t^4+8t^3+5t^2+t-1)=0$ which clearly has no integer solutions.
\end{case}

Thus we see that $\Th_2$ cannot be real.

\begin{case}
 Now suppose $s=4$.  Then we obtain a polynomial with irreducible factorization: $$(2t+1)(4t^2-4t-1)(4t^6+12t^5+7t^4-6t^3-11t^2-6t-1)$$
which clearly has no integer roots (or note that it is an odd number for any integer $t$).
\end{case}
\begin{case}
If $s=3$ then $t,v$ must satisfy:
$$(1+7t+16t^2+18t^3+9t^4)(4t^2+2t+1-8v^2),$$
for which both factors are odd for integer $t,v$.
\end{case}
\begin{case}
 If $s=6$ then $t$ must satisfy:
$$(8t^2-4t-1)(128t^8+512t^7+904t^6+920t^5+573t^4+210t^3+36t^2-t-1)$$
which again we see is an odd number for any integer $t$.
\end{case}

\begin{case}
 If $s=7$ then $t$ must satisfy 
$$t(1+t)(8t^3+4t^2-4t-1)(8t^3-12t^2-8t-1)P_{18}(t)$$
where $P_{18}(t)$ is a polynomial of degree $18$ that takes odd values.
We conclude that in this case $t=0$ is the only non-negative integer solution.  We have seen that $t=0$ implies $\Th_2$ is a $7$th root of unity, so this uniquely determines all $n_i$, $d,f$ and $g$ as well as $\Th_1$ and $\Th_3$ up to the particular choice of a $7$th root of unity $\Th_2$.
\end{case}

\begin{case}
 If $s=9$ then $t$ must satisfy the following polynomial modulo $2$:
$$(t^6-t^3+1)(t^6+t^4+t^3+t+1),$$
which is clearly odd.
\end{case}

\begin{case}
 If $s=14$ then $t$ must satisfy:
$$(8t^3+4t^2-4t-1)(64t^6+192t^5-208t^4-64t^3+32t^2+12t+1)Q_{24}(t)$$
where $Q_{24}(t)=2^{12}t^{24}+\cdots+2^3$.  The first two factors are clearly odd and so cannot be zero.  Any non-negative integer roots of the polynomial $Q_{24}(t)$ must be $1,2,4$ or $8$ by the rational-root theorem, each of which can be eliminated. 
\end{case}

\begin{case}
 If $s=18$ then $t$ must satisfy:
$$(24t^3-6t-1)(64t^6-192t^4-32t^3+36t^2+12t+1)R_{24}(t)$$
where $R_{24}(t)$ is a polynomial of degree $24$.  Each factor takes odd values
for integral $t$, so there are no solutions.
\end{case}

This completes the proof of Theorem \ref{Z6theorem}.
\begin{theorem}
The category with $\Gal(\CC)$ isomorphic to $\lan (0\; 1), (3\; 4)\ra$
does not exist.
\end{theorem}

\begin{proof}
Consider the $\tS$-matrix in eqn. (\ref{s-mat}) and  
let $\sigma=(0\; 1)\in S_5$. From the
 formula
 $\tilde{s}_{j,k}=\epsilon_{\tau(j),\tau}\epsilon_{k,\tau}\tilde{s}_{\tau(j),\tau^{-1}(k)}$
 for each $\tau\in \Gal(\CC)$ in Proposition \ref{rswprop}, we have
 the following:

 \begin{eqnarray}
&&1=\tilde{s}_{0,0}=\epsilon_{1,\sigma}\epsilon_{0,\sigma}x\label{ep1}\\
&&b=\tilde{s}_{3,2}=\epsilon_{3,\sigma}\epsilon_{2,\sigma}b\label{ep2}\\
&&f=\tilde{s}_{2,0}=\epsilon_{2,\sigma}\epsilon_{0,\sigma}y\label{ep3}\\
&&y=\tilde{s}_{2,1}=\epsilon_{2,\sigma}\epsilon_{1,\sigma}f\label{ep4}\\
&&g=\tilde{s}_{3,0}=\epsilon_{3,\sigma}\epsilon_{0,\sigma}a\label{ep5}
\end{eqnarray}

\noindent Then we notice two facts. One is
 $\epsilon_{0,\sigma}=\epsilon_{1,\sigma}$ from (\ref{ep3}) and (\ref{ep4}), and thus
$x=1$ from (\ref{ep1}). The other is: provided $b\not=0$,
$\epsilon_{2,\sigma}=\epsilon_{3,\sigma}$ from (\ref{ep2}), thus
$\epsilon_{2,\sigma}\epsilon_{0,\sigma}=\epsilon_{3,\sigma}\epsilon_{0,\sigma}$,
which means we have two possibilities: $y=f$ and $a=g$, or $y=-f$
and $a=-g$ from (\ref{ep3}) and (\ref{ep5}).  If $b=0$, (\ref{ep2}) gives no information, and the analysis is more delicate.

We first assume that $b\not=0$.

For the case that $y=f$ and $a=g$, we have an easy contradiction
as follows:
$\tS^2=D^2C$ implies that
$2d+f^2+2g^2=0$ which contradicts $d>0$.

Now we consider the other case that $y=-f$ and $a=-g$.
The  gives us the following:
\begin{eqnarray}
 &&D^2=1+d^2+f^2+2g^2\label{or1}\\
&&2d=f^2+2g^2\label{or2}\\
&&f(d-1-z)=2gb\label{or3}\\
&&D^2=z^2+2(b^2+f^2)\label{or4}\\
&&2(g^2+h_1^2-h_2^2)+b^2=0\label{or5}\\
&&g(1-d+2h_2)+fb=0\label{or6}\\
&&2fg+zb+2bh_1=0\label{or7}\\
&&2(g^2+h_1^2+h_2^2)+b^2=D^2\label{or8}
\end{eqnarray}

From (\ref{or1}) and (\ref{or2}) we obtain $D=d+1$ (assuming $D>0$ and $d\geq 1$).
By replacing $\tS$ by its complex conjugate if necessary, we may assume that $h_2>0$ so that (\ref{or5}) and (\ref{or8}) imply $h_2=D/2$.  To derive further consequences of these equations we use Maple's Gr\"obner basis algorithm with a monomial order that eliminates $D$ and $g$. 
One consequence is:
\begin{equation}
 (2h_1(D-1)+2g^2+D-D^2)(2h_1-2g^2+D)=0
\end{equation}
so we have two possibilities:
$h_1=\frac{-2g^2-D+D^2}{2D-2}$, or $h_1=g^2-\frac{D}{2}$.

For the first case, we obtain as a consequence $(1-D+g^2)(fg+Db-b)=0$. The first
factor is non-zero since (\ref{or2})
implies $D-1=d=g^2+\frac{1}{2}f^2$ so $D>1+g^2$. Thus we must have $fg=b(1-D)=-bd$. Notice that
$\frac{b}{f}=-\frac{g}{d}$ are both integers so that
$-\frac{g}{d}=\sigma(-\frac{g}{d})=g$. We have a contradiction
$d=-1$.

For the other case, we obtain $(1-D+g^2)(fg-b)=0$ which implies
$fg=b$ from the same reasoning as above. This means that
$g=\frac{b}{f}$ is an integer. Since $\sigma$ permutes the roots
$g$ and $-\frac{g}{d}$ we get $d=-1$.

We are left with two remaining cases: $b=0$ and (Case 1) $y=f$ and $a=-g$ or  (Case 2) $y=-f$ and $a=g$.  Although these cases can likely be eliminated without resorting to symbolic computation techniques, we follow the strategy used in the proof of Theorem \ref{Z6theorem} to further illustrate its efficacy. We first consider Case 1 as it is slightly more involved. 
The $S$-matrix is:
$$ \tS=\begin{pmatrix}
          1 & d & f & g & g\\
d & 1 & f& -g & -g\\
f & f & z& 0 & 0\\
g & -g & 0& h & \overline{h}\\
g & -g & 0& \overline{h} & h\end{pmatrix} $$

  Clearly $g$ cannot be an integer, since $Gal(\CC)$ interchanges $g$ and $-g/d$.  similarly, $f\in\N$ if and only if $d\in\N$ if and only if $d=1$.  The equation $\tS^2=D^2C$ immediately implies:
\begin{eqnarray}
 &&z=-(1+d)\\
&&h_1=\frac{d-1}{2}\\
&&d=g^2-f^2/2 \label{b=0orthrel1}
\end{eqnarray}
Notice that if $d=1$ then $z=-2$ and $z/f\in\Z$ implies $f=1$ or $f=2$.  But eqn. (\ref{b=0orthrel1}) and the fact that $g^2$ is an algebraic integer imply that if $d$ or $f$ are integral, then $d=1$, $f=2$, $g=\sqrt{3}$.  We will consider this case seperately; for now assume that $d$ is not an integer.

Next we consider the Diophantine equations satisfied by the entries $n_i$ of the fusion matrices $N_1$, $N_2$ and $N_3$ from (\ref{N1}-\ref{N3}).  In this case we have a common integer eigenvalue for the $N_i$, namely $\psi_2=(1,1,-m,0,0)$ where $m=(d+1)/f$.  Besides the commutation relations, we have relations obtained from the characteristic polynomials and the relations involving $m$.
We also obtain $$(n_{10}-n_{11})(m+n_{10}-n_{11})=(n_6-n_7)^2-1=0$$ 
from the fact that the only integer eigenvalues of $N_1$ are $\pm 1$ and the only integer eigenvalues of $N_2$ are $0$ and $-m$ (provided $d\not=1$).  From these relations we can solve for all of the $N_i$ entries and $m$ in terms of four variables $u:=n_8$, $v:=n_{11}$, $t:=n_{13}$ and $w:=n_{14}$ which satisfy the three relations:
\begin{eqnarray*}
 &&(t-w)(t^3+t^2w+2t-2tv^2-tw^2+2wv^2+2w-w^3),\\ &&(t-w)(-3vw-vt+ut+uw),\\
 &&2+t^2-w^2+2uv-4v^2.
\end{eqnarray*}
Unfortunately these Diophantine equations appear to have many solutions.
We proceed to analyze the relations implied by eqns. (\ref{threls}):
\begin{eqnarray}
 &&\theta_1^2-1-n_1d\theta_1-n_2f\theta_2-2n_3g\theta_3,\\ &&\theta_1\theta_2f-n_2d\theta_1-n_4f\theta_2-2n_5g\theta_3,\\ &&\theta_1\theta_3g+n_3d\theta_1+n_5f\theta_2+(n_6+n_7)g\theta_3,\\ &&n_5d\theta_1+n_9f\theta_2+(n_{10}+n_{11})g\theta_3,\label{b0threl4}\\ &&\theta_2^2(d+1)+1+n_4d\theta_1+n_8f\theta_2+2n_9g\theta_3,\\ &&\theta_3^2(d-1)-1-d\theta_1(n_7+n_6)-f\theta_2(n_{11}+n_{10})-g\theta_3(2n_{12}+n_{13}+n_{14})
\end{eqnarray}
as well as the relations implied by the common positive eigenvector $\psi_0$.
By combining the relations involving $\theta_i$ as well as the $S$-matrix and $N_i$ entries we immediately obtain $\theta_1=1$.  From this we are able to eliminate $\theta_3$ (assuming $n_5\not=0$, which implies $d=1$) and find a polynomial of degree two in $\theta_2$ with coefficients in $\Q(d)=\Q(d,f,g)$.  Thus $[\Q(\theta_2):\Q]\in\{1,2,4\}$ so that $\theta_2$ has degree $1\leq k\leq 6$ or $k\in\{8,10,12\}$.  Using a Gr\"obner basis algorithm with pure-lexicographic order $d<f<g<\cdots<t<w$ on the full set of relations together with the appropriate cyclotomic polynomial in $\theta_2$, we eliminate each of these cases in turn.  In each $k\geq 3$ we obtain as a consequence a product of irreducible polynomials of degree at least $2$ in $w$ with the Diophantine equation $(t-w)^2+4$, which has no solutions.  In case $k=1,2$ we obtain $d=1$.

So as soon as we rule out $d=1$, we are done with Case 1.  We find that in this case we have:
\begin{eqnarray*}
 &&n_1= n_2 = n_3  = n_5 =n_6  = n_9 = n_{12} = n_{13} = n_{14} = 0,\\ &&n_4 =n_7 =n_8 =n_{10}=n_{11}= 1.
\end{eqnarray*}
From (\ref{b0threl4}) we immediately obtain $\theta_3=0$, a contradiction.

Case 2 is similar, but slightly easier.  Here we have $S$-matrix:
$$\tS=\begin{pmatrix}
          1 & d & f & g & g\\
d & 1 & -f& g & g\\
f & -f & z& 0 & 0\\
g & g & 0& h & \overline{h}\\
g & g & 0& \overline{h} & h\end{pmatrix}.$$
Since $\Gal(\CC)$ interchanges $\phi_0$ and $\phi_1$ we see that in this case if any of $g,f$ or $d$ are integral, $f=0$, a contradiction.  The rest of the proof is essentially the same as for Case 1: one finds that $\Th_1=1$ and $[\Q(\Th_2):\Q]$ must be $1,2$ or $4$ polynomial over $\Q$.  Going through the cases one finds that there are no integer solutions.
\end{proof}

\section{Conclusions and Future Directions}\label{conclusions}

We briefly describe how our computational strategy can potentially be employed to classify modular categories of fixed (low) rank $n$, or at least verify Conjecture \ref{wangconj} for $n$.  First one classifies the integral modular categories, perhaps using Lemma \ref{egyptlemma}.  Next one determines all abelian subgroups of $S_n$ that do not fix the label $0$ and considers them case-by-case.  From Proposition \ref{rswprop} and the equation $\tS^2=D^2C$ one determines the form of the $\tS$-matrix, eliminating as many of the $\tS$-matrix variables as possible.  Many cases will be eliminated in this way.  For each remaining case, one uses the fusion matrices and the coefficients of their characteristic polynomials to obtain Diophantine equations and eliminate integer variables using a Gr\"obner basis algorithm.  From relation (\ref{threls}) one determines the possible degrees of the roots of unity $\theta_i$.  Finally one uses the relations among 
the $\tS$-matrix entries, the fusion multiplicities and the $\theta_i$ together with assumptions on the degree of the $\theta_i$ to obtain new Diophantine equations.  These new Diophantine equations will often give finitely many solutions.

In this work we have classified (the Grothendieck semirings of) modular categories under the restrictions: 1) some object is non-self-dual, 2) the dimensions of simple objects are positive and 3) rank $\leq 5$.  The strategy described above does not require any of these restrictions, and we make some remarks on the complexity of removing these restrictions.

\begin{enumerate}
 \item To classify rank $5$ non-integral pseudo-unitary self-dual modular categories one must consider around $9$ Galois groups, only two of which are known to be realized for a modular category.
\item The pseudo-unitarity assumption can probably be removed without loss of generality for rank $\leq 5$, although this assumption is often physically justified.  For rank $6$, it is known that there are modular fusion algebras that have no pseudo-unitary realization (see  \cite{Rsurvey}) but are realized as $Gr(\CC)$ for a modular category $\CC$ nonetheless.
\item In rank $6$ one encounters product categories, and at least $9$ modular fusion algebras are known to exist. The (two) known non-self-dual examples are products of the unique rank $3$ non-self-dual fusion algebra with each of the two rank $2$ fusion algebras.  Under assumptions 1) and 2) the classification could likely be carried out to rank $7$ using our methods.
\end{enumerate}

Finally, while the assumption of modularity seems critical to our approach, analogous results can presumably be obtained for \emph{pre-modular} categories, in which one omits the invertibility of the $S$-matrix.  Brugui\`eres \cite{Br} has shown that if a pre-modular category fails to be modular then two columns of the $S$-matrix are proportional so that Galois group identities might still be obtained.  Moreover, eqn. (\ref{threls}) and various other identities hold generally in a pre-modular category so that a computational approach is still viable.

\end{document}